\documentclass[onecolumn]{article}
\usepackage[textwidth=14cm]{geometry}

\usepackage{listings}

\usepackage{graphicx}
\usepackage{tikz}

\usepackage{amssymb,amsmath}
\usepackage{enumitem}      
\usepackage{url}

\newcommand{\map}[3]{#1: #2 \rightarrow #3}
\renewcommand{\natural}{{\mathbb{N}}}

\newcommand{\until}[1]{\{1,\ldots,#1\}}

\newcommand{\st}{\text{subject to }}

\newcommand{\m}{\mathop{\textrm{minimize}}}

\newcommand{\R}{\mathbb{R}}

\def\disropt/{\textsc{disropt}}

\definecolor{mygray}{gray}{0.95}
\graphicspath{{figs/}}

\title{DISROPT: a Python Framework \\for Distributed Optimization\footnote{This result is part of
  a project that has received funding from the European Research Council (ERC) under
  the European Union's Horizon 2020 research and innovation programme
  (grant agreement No 638992 - OPT4SMART).}}

\author{Francesco Farina, Andrea Camisa, Andrea Testa, \\ Ivano Notarnicola, Giuseppe Notarstefano}

\date{\small Department of Electrical, Electronic and Information Engineering,\\
Alma Mater Studiorum Universit\`{a} di Bologna, Bologna, Italy\\
  $\{$\texttt{franc.farina}, \texttt{a.camisa}, \texttt{a.testa},
  \texttt{ivano.notarnicola}, \texttt{giuseppe.notarstefano}$\}$\texttt{@unibo.it}
}

\begin{document}

\maketitle

\begin{abstract}
    In this paper we introduce \disropt/, a Python package for distributed optimization
    over networks. We focus on cooperative set-ups in which an optimization problem
    must be solved by peer-to-peer processors (without central coordinators) that have
    access only to partial knowledge of the entire problem.
    To reflect this, agents in \disropt/ are modeled as entities that are initialized with their
    local knowledge of the problem. Agents then run local routines and communicate
    with each other to solve the global optimization problem.
    A simple syntax has been designed to allow for an easy modeling of the problems.
    The package comes with many distributed optimization algorithms
    that are already embedded. Moreover, the package provides full-fledged
    functionalities for communication and local computation, which can be
    used to design and implement new algorithms.
    \disropt/ is available at \url{github.com/disropt/disropt} under the GPL license,
    with a complete documentation and many examples.
\end{abstract}


\section{Introduction}

In recent years, distributed learning and control over networks have gained a growing attention. Most problems arising in these contexts can be formulated as distributed optimization problems, whose solution calls for the design of tailored strategies.
The main idea of distributed optimization is to solve an optimization problem
over a network of computing units, also called agents.
Each agent can perform local computation and can exchange information only with its neighbors in the network.
Typically, the problem to be solved is assumed to have a given structure, while the
communication network can be unstructured. Each agent knows only a portion of
the global optimization problem, so that finding a solution to the network-wide
problem requires cooperation with the other agents.
A distributed algorithm consists of an iterative procedure in which each agent
alternates communication and computation phases with the aim of eventually
finding a solution to the problem.
The recent monograph~\cite{notarstefano2019distributed} provides a comprehensive
overview of the most common approaches for distributed optimization, together with
the theoretical analysis of the main schemes in their basic version.

In this paper we present \disropt/, a Python package designed to run distributed
optimization algorithms over peer-to-peer networks of processors.
%
%
In the last years, several toolboxes have been developed in order to solve optimization
problems using \emph{centralized} algorithms.
Examples of toolboxes written in C are \textsc{osqp}~\cite{stellato2018osqp}, 
and \textsc{GLPK}~\cite{makhorin2008glpk}.
As for packages developed in C\texttt{++}, nonlinear optimization problems can be solved
by using \textsc{opt\texttt{++}}~\cite{meza1994opt++}.
In the context of nonlinear optimization, we mention
\textsc{acado}~\cite{Houska2011a}, which deals with optimal control,
 and \textsc{ipopt}~\cite{biegler2009large}, which solves large-scale problems.
More recent interpreted languages, such as Matlab and Python, do not have the
performance of low-level compiled languages such as C and C\texttt{++}, however
they are more expressive and often easier to use.
Well-known Matlab packages for optimization are \textsc{yalmip}~\cite{lofberg2004yalmip}
and \textsc{cvx}~\cite{cvx}.
Notable Python packages for convex optimization are \textsc{cvxpy}~\cite{cvxpy}
and \textsc{cvxopt}~\cite{andersen2015cvxopt}.
Based on \textsc{cvxpy}, the toolbox \textsc{snapvx}~\cite{hallac2017snapvx} allows for
the solution of large-scale convex problems defined over graphs by exploiting their structure.
An extension of~\cite{cvxpy} to optimize convex objectives over nonconvex domains using
heuristics is \textsc{ncvx}~\cite{diamond2018general}.
Other well-known packages based on the recent programming language Julia are \textsc{optim},
\textsc{convex.jl} and \textsc{jump}~\cite{mogensen2018optim,udell2014convex,dunning2017jump}.
A Julia package for stochastic optimization is~\cite{huchette2014parallel},
which implements a parallel solver.

None of the above references provides direct capabilities to solve optimization problems
over networks using distributed algorithms.
The aim of \disropt/ is to bridge this gap.
The package is designed as follows. When a distributed algorithm is executed,
many (possibly spatially distributed) processes are created.
Each process corresponds to an agent in the network, has its own memory space
with its private data, runs its own set of instructions and cooperates with
the other agents through a message-passing paradigm.
Consistently with the distributed model, there is no central coordinator.
The package provides a comprehensive framework
for all the typical tasks that must be performed by distributed optimization algorithms.
In particular, the package allows for both synchronous and asynchronous communication over
custom networks of peer-to-peer agents.
Moreover, it provides an easy-to-use interface to represent, for each agent, the local knowledge of
the global optimization problem to be solved and to run distributed optimization algorithms
via a streamlined interface.
Local optimization problems can be also solved (if needed).
The tools provided by \disropt/ let the user to easily design new
distributed optimization algorithms by using an intuitive syntax and to solve several classes of problems arising both in distributed control and machine learning frameworks.

The paper is organized as follows. In Section~\ref{sec:setup}, the distributed
optimization framework is introduced. The architecture of \disropt/ is presented in
Section~\ref{sec:architecture} and a numerical computation on three example scenarios
is provided in Section~\ref{sec:example}.



\section{Distributed Optimization Set-ups}
\label{sec:setup}

In this section, we introduce the optimization set-ups considered in \disropt/.
Formally, an optimization problem is a mathematical problem which consists
in finding a minimum of a function while satisfying a given set of constraints.
In symbols,
\begin{equation*}
    \begin{aligned}
        &\m_{x}
        & & f(x) \\
        & \st
        & & x \in X,
    \end{aligned}
\end{equation*}
where $x \in \R^d$ is called optimization variable, $\map{f}{\R^d}{\R}$
is called cost function and $X \subseteq \R^d$ describes the problem constraints.
The optimization problem is assumed to be feasible, to have finite optimal cost and to admit at
least an optimal solution, which is usually denoted as $x^\star$. The optimal solution is a
vector satisfying all the constraints and attaining the optimal cost. If the problem is nonconvex, $x^\star$ can be any (feasible) stationary point.

Distributed optimization problems arising in applications usually enjoy a proper
structure in their mathematical formulation. In \disropt/, three different optimization
set-ups are considered and are detailed next.

\subsection{Cost-coupled Set-up}
\label{sec:cost-coupled}
In the \emph{cost-coupled optimization set-up}, the cost function is expressed as the sum of local contributions $f_i$
and all of them depend on a common optimization variable $x$. Formally, the set-up is
\begin{equation}
    \begin{aligned}
        &\m_{x}
        & & \sum_{i=1}^N f_i(x) \\
        & \st
        & & x \in X,
    \end{aligned}
\label{eq:pb_cost_coupled}
\end{equation}
where $x \in \R^d$ and $X \subseteq \R^d$. The global constraint set $X$ is
known to all the agents, while $\map{f_i}{\R^d}{\R}$ is assumed to be known
by agent $i$ only, for all $i\in \until{N}$.

In some applications, the constraint set $X$ can be expressed as the intersection
of local constraint sets, i.e.,
\begin{equation*}
    X = \bigcap_{i=1}^N X_i,
\end{equation*}
where each $X_i \subseteq \R^d$ is assumed to be known by agent $i$ only,
for all $i\in \until{N}$.

The goal for distributed algorithms for the cost-coupled set-up is that all the
agent estimates of the optimal solution of the problem are eventually consensual to $x^\star$.


\subsection{Common-cost Set-up}\label{sec:common-cost}
In the \emph{common-cost optimization set-up}, there is a unique cost function $f$ that depends on a common
optimization variable $x$, and the optimization variable must further satisfy local constraints.
Formally, the set-up is
\begin{equation}\label{pb:common_cost}
    \begin{aligned}
        &\m_{x}
        & & f(x)\\
        & \st
        & & x \in \bigcap_{i=1}^N X_i
    \end{aligned}
\end{equation}
where $x \in \R^d$ and each $X_i \subseteq \R^d$. The cost function $f$
is assumed to be known by all the agents, while each set $X_i$ is known by agent
$i$ only, for all $i\in\until{N}$.

The goal for distributed algorithms for the common-cost set-up is that all the
agent estimates of the optimal solution of the problem are eventually consensual to $x^\star$.


\subsection{Constraint-coupled Set-up}\label{sec:constraint-coupled}

In the \emph{constraint-coupled optimization set-up}, the cost function is
expressed as the sum of local contributions $f_i$ that depend on a local
optimization variable $x_i$. The variables must satisfy local constraints
and global coupling constraints among all of them. Formally, the set-up is
\begin{equation}
    \begin{aligned}
        &\m_{x_1,\ldots,x_N}
        & & \sum_{i=1}^N f_i(x_i)\\
        & \st
        & & x_i \in X_i, \quad i \in \{1, \ldots, N\}\\
        & & & \sum_{i=1}^N g_i(x_i) \leq 0
    \end{aligned}
\end{equation}
where each $x_i \in \R^{d_i}$, $X_i \subseteq \R^{d_i}$, $\map{f_i}{\R^{d_i}}{\R}$
and $\map{g_i}{\R^{d_i}}{\R^S}$, for all $i \in \until{N}$.
Here the symbol $\leq$ is also used to denote component-wise inequality for vectors.
Therefore, the optimization variable consists of the stack of all the variables $x_i$,
namely the vector $(x_1,\ldots,x_N)$.
The quantities with the subscript $i$ are assumed to be known by agent $i$ only,
for all $i\in\until{N}$. The function $g_i$, with values in $\R^S$, is used to express
the $i$-th contribution to $S$ coupling constraints among all the variables.

The goal for distributed algorithms for the constraint-coupled set-up is that each agent asymptotically
computes its portion $x_i^\star \in X_i$ of an optimal solution $(x_1^\star, \ldots, x_N^\star)$
of the optimization problem, thereby satisfying also the coupling constraints
$\sum_{i=1}^N g_i(x_i^\star) \leq 0$.

\section{Software Architecture and Basic Syntax}
\label{sec:architecture}
\begin{figure}[htbp]\centering
  \includegraphics[scale=0.9]{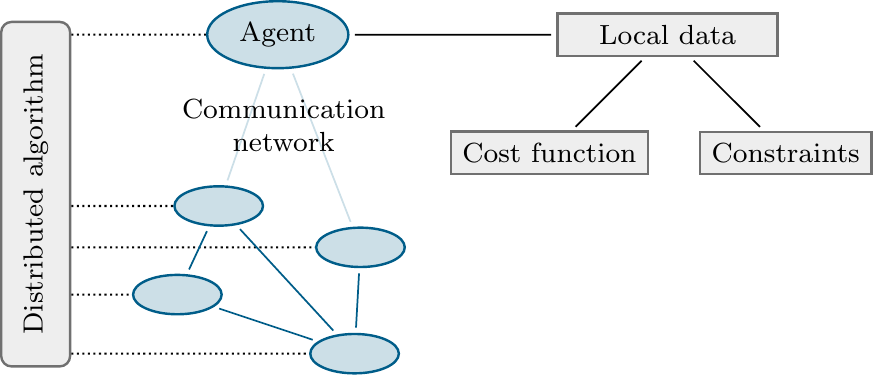}
  \caption{
    Distributed scenario architecture.
    Agents are equipped with their local information and interact with
    each other through the communication network to run the
    distributed algorithm.
  }
  \label{fig:architecture}
\end{figure}

In this section, we present the architecture of the package, which
replicates the typical structure of a distributed algorithm, see
Figure~\ref{fig:architecture}. The main entities of a distributed scenario
are the agents (with their local information), the communication network,
and the local routines of the distributed algorithm.

In \disropt/, these components are represented with an object-oriented framework.
In the remainder of this section, we provide a brief description of the main
classes and of the semantics of \disropt/.

\subsection{Agent}
An instance of the \emph{Agent} class represents a single computing unit in the network.
This class is equipped with the list of the neighboring agents.
It is also embedded with an instance of the class \emph{Problem} (detailed next),
which describes the local knowledge of the global optimization problem.

\begin{figure}[!htbp]
	\centering
	\includegraphics{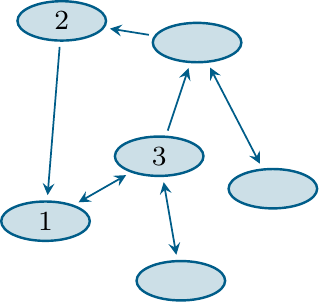}
	\caption{Example of a network of $6$ agents.}
	\label{fig:network_ex}
\end{figure}

Suppose that we want to instantiate the agent with index $1$ of the network in Figure~\ref{fig:network_ex}.
The in-neighbors of agent $1$ are agents $2$ and $3$ while the unique out-neighbor is agent $3$.
The following Python code can be used to accomplish the task:
\begin{lstlisting}
from disropt.agents import Agent

in_nbrs = [2, 3]
out_nbrs = [3]
agent = Agent(in_nbrs, out_nbrs)
\end{lstlisting}

\subsection{Communication Network}
The \emph{Communicator} class handles communication among the network agents.
It allows agents to send and receive data from neighbors in a synchronous/asynchronous,
time-invariant/time-varying fashion.
In the current release, communication over the Message Passing Interface (MPI)
is supported, however custom communication protocols can be implemented as well.

As an example, to perform synchronous communication to exchange
a two-dimensional vector with neighbors, the syntax is
\begin{lstlisting}
from disropt.communicators import MPICommunicator

vect = numpy.random.rand(2, 1)
comm = MPICommunicator()
exch_data = comm.neighbors_exchange(vect, in_nbrs, out_nbrs, dict_neigh=False)
\end{lstlisting}
The flag \verb|dict_neigh| set to \verb|False| means that the same object
is sent to all the neighbors, while \verb|exch_data|
is a dictionary containing all the vectors received
from the neighbors with their corresponding indices.

An instance of this class is embedded in the Agent class which also manages
the list of in- and out-neighbors. Therefore, the previous code can also be
restated as follows:
\begin{lstlisting}
vect = numpy.random.rand(2, 1)
exch_data = agent.neighbors_exchange(vect)
\end{lstlisting}

Note that the provided communication features are typically required
only during the algorithm implementation.

\subsection{Distributed Algorithm}
The class \emph{Algorithm} aims to characterize the behavior and the local steps of
distributed algorithms. This class reads the local problem data, handles the data
received by neighbors and updates the local solution estimate.
Specializations of this class correspond to different distributed algorithms.
We implemented several algorithms corresponding to the different distributed
optimization set-ups of Section~\ref{sec:setup}.
An exhaustive list of the currently implemented
algorithms is provided in Table~\ref{tb:algos}. References to the implemented algorithms can be found in~\cite{notarstefano2019distributed}.

\begin{table*}[!t]
\caption{Overview of implemented algorithms. $^*$=next release.}
\vspace{2ex}
\scalebox{0.63}{
\begin{tabular}{|l|c|c|c|c|c|c|c|}
\hline
\bf Algorithm           & \bf Communication  & \bf Time-varying   & \bf Block-wise        & \bf Constraints   & \bf Non-smooth & \bf Asynchronous  \\
\hline\hline
\multicolumn{7}{l}{\bf Cost-coupled } \\
\hline\hline
Distributed subgradient          & Directed  & \checkmark   &  \checkmark       & Global        & \checkmark & \checkmark \\ \hline
Gradient tracking               & Directed       & \checkmark &  \checkmark $^*$  & No       & &  \\ \hline
Distributed dual decomposition        & Undirected   &         &                   & Local   & \checkmark & \checkmark$^*$  \\ \hline
Distributed ADMM               & Undirected   &         &                   & Local   & \checkmark &  \\ \hline
ASYMM                     & Undirected   &         &  \checkmark       & Local   & & \checkmark \\ \hline
\multicolumn{7}{l}{\bf Common cost }\\
\hline\hline
Constraints consensus              & Directed       & \checkmark &                   & Local   & & \checkmark $^*$  \\ \hline
Distributed set membership             & Directed       & \checkmark   &                   & Local   & & \checkmark  \\ \hline
\multicolumn{7}{l}{\bf Constraint-coupled }\\
\hline\hline
Distributed Dual Subgradient   & Directed  &         &                   & Local   & \checkmark &  \\ \hline
Distributed Primal Decomposition  & Undirected   &         &                   & Local   & \checkmark &  \\ \hline
\end{tabular}
}

\label{tb:algos}
\end{table*}

For example, in order to run the distributed subgradient method
for $100$ iterations,
the Python code is:
\begin{lstlisting}
from disropt.algorithms import SubgradientMethod

x0 = numpy.random.randn(2, 1)
algorithm = SubgradientMethod(agent=agent, initial_condition=x0)
algorithm.run(iterations=100)
\end{lstlisting}
Notice that we are assuming that the \verb|agent| is already equipped with the
local problem information, which can be done as described in the next subsection.

Then, it is possible to run the distributed algorithm over a network of
$N$ agents by executing:
\begin{lstlisting}
mpirun -np N python <source_file.py>
\end{lstlisting}
%

\subsection{Local Agent Data}
%
%
%
%

The class \emph{Problem} is used to model the locally available data of the
global optimization problem.
It is embedded with an objective function
(an instance of the class \emph{Function}) and a list of constraints
(a list of objects of the class \emph{Constraint}).
The class should be initialized according to the specific distributed
optimization set-up and must be provided to the class \emph{Agent}.

For instance, suppose that, in a cost-coupled set-up, the agent knows the following
function and constraint,
\begin{align*}
	f_i(x) = \|x\|^2,
	\hspace{1cm}
	X_i = \{x \in \mathbb{R}^2 \mid -1 \le x \le 1\}.
\end{align*}

The corresponding Python code is:
\begin{lstlisting}
from disropt.functions import SquaredNorm, Variable
from disropt.problems import Problem

x = Variable(2)
objective_function = SquaredNorm(x)
constraints = [x >= -1, x <= 1]
problem = Problem(objective_function, constraints)
agent.set_problem(problem)
\end{lstlisting}

Since many distributed algorithms require the solution
of small optimization problems, this class is also able to solve the
problem described by the cost function and the constraints
through the method \verb|solve()|, which currently relies upon \textsc{cvxopt},
\textsc{glpk}, \textsc{ospq} and \textsc{cvxpy}.

\section{Example Usage}\label{sec:example}
In this section we present three illustrative examples on how \disropt/ can be used to
solve a cost-coupled (cf. Sec.~\ref{sec:cost-coupled}), a common-cost (cf. Sec.~\ref{sec:common-cost})
and a constraint-coupled optimization problem (see Sec.~\ref{sec:constraint-coupled}), respectively.
In particular, we consider a distributed classification problem where the training
points are spread among the agents and a distributed microgrid control problem.

\subsection{Distributed Classification via Logistic Loss}\label{sec:log_loss}
\begin{figure}[htpb]
\centering
  \includegraphics[scale=1]{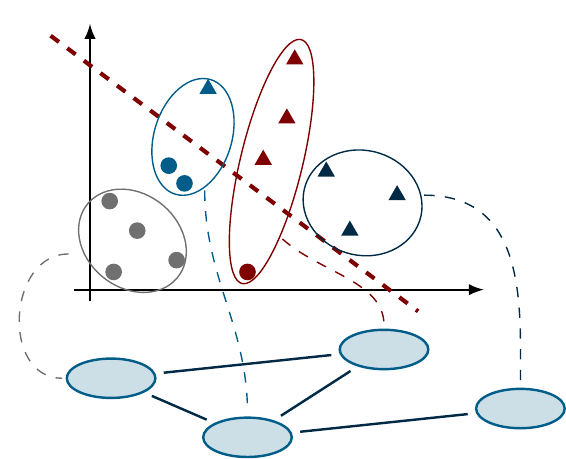}
  \caption{Illustration of a distributed classification problem.}
\label{fig:classification}
\end{figure}

The classification problem consists in dividing a set
of points, representing data in a feature space, into two clusters, by means of
a separating hyperplane. The purpose of the agents is to cooperatively
estimate the parameters of the hyperplane.
In Figure~\ref{fig:classification}, a bidimensional example is reported. The classification
task corresponds to computing the parameters of the line (in red) separating triangles from circles.

\subsubsection{Problem formulation}

Let us consider $N$ agents contributing to the $d$-dimensional classification problem as follows.
Each agent $i$ is equipped with $m_i$ points $p_{i,1}, \ldots, p_{i,m_i} \in \R^d$. 
The points are associated to binary labels, that is, each point $p_{i,j}$ is labeled with
$\ell_{i,j} \in \{-1,1\}$, for all $j \in \{1, \ldots, m_i\}$ and $i \in \{1, \ldots, N\}$.
The problem consists of building a linear classification model from the training
samples by maximizing the a-posteriori probability of each class.
In particular, we look for a separating hyperplane of the form
$\{ z \in \mathbb{R}^d \mid w^\top z + b = 0 \}$, whose parameters
$w$ and $b$ can be determined by solving the following optimization problem
\begin{equation*}
    \begin{aligned}
        &\m_{w, b}
        & & \sum_{i=1}^N \: \sum_{j=1}^{m_i}
        \log \left( 1 + e^{ -(w^\top p_{i,j} + b) \ell_{i,j} } \right) + \frac{C}{2} \|w\|^2,
    \end{aligned}
\end{equation*}
where $C > 0$ is the regularization parameter. As we already mentioned, this
is a cost-coupled problem of the type~\eqref{eq:pb_cost_coupled},
where each local cost function $f_i$ is appropriately defined\footnote{The regularization term can be appropriately split among the agents so that $f_i(x) = \sum_{j=1}^{m_i} \log \left( 1 + e^{ -(w^\top p_{i,j} + b) \ell_{i,j} } \right) + \frac{C}{2N} \|w\|^2$.}.
%
The goal is to make agents agree on a common solution $(w^\star, b^\star)$,
so that all of them can compute the separating hyperplane as
$\{ z \in \mathbb{R}^d \mid (w^\star)^\top z + b^\star = 0 \}$.

%


\subsubsection{Simulation results}

We consider a bidimensional sample space ($d = 2$).
Each agent $i$ generates a total of $m_i$ points of both labels, with $m_i$ a random
number between $4$ and $10$.
For each label, the samples are drawn according to a multivariate Gaussian distribution,
with covariance matrix equal to the identity and mean equal to $(0,0)$
for the label $1$ and $(3,2)$ for the label $-1$.
The regularization parameter is $C = 10$.


We run a comparative study of the distributed subgradient algorithm and
the gradient tracking algorithm, with $N = 20$ agents and $20\,000$ iterations.


As for the step size, we use the following rules: constant step-size $\alpha^t = 0.001$
for gradient tracking and diminishing step-size $\alpha^t = ( 1/t )^{0.6}$ 
or distributed subgradient.

%

\begin{figure}[htpb]\centering
  \includegraphics[scale=1]{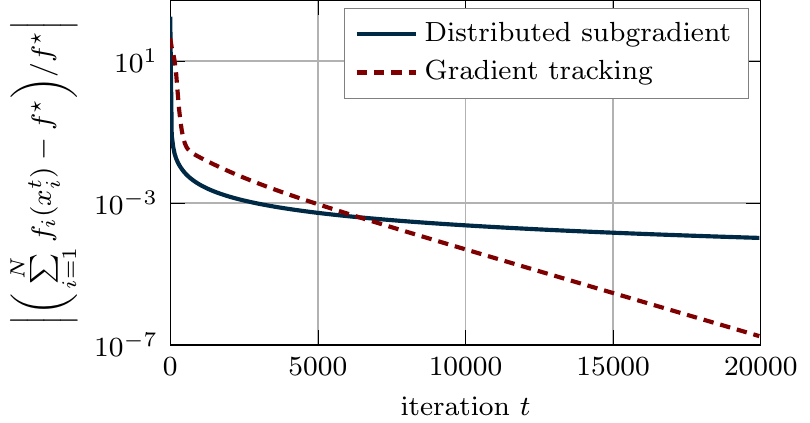}
  \caption{Distributed classification: normalized cost error between the locally computed
    solution estimates and the optimal cost.}
  \label{fig:logistic_cost}
\end{figure}

\begin{figure}[htpb]\centering
  \includegraphics[scale=1]{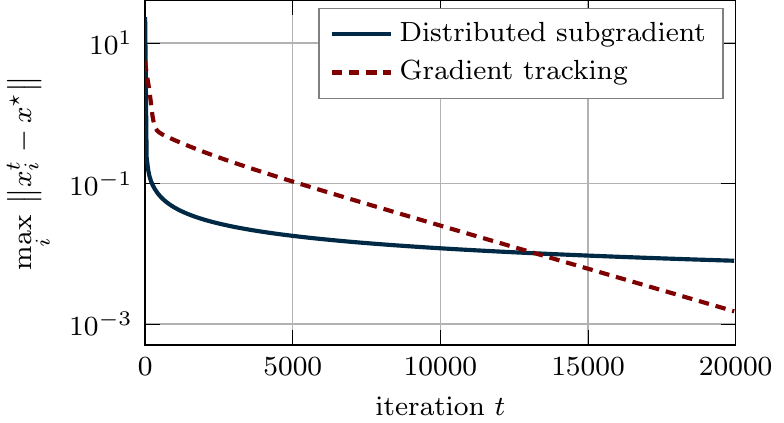}
  \caption{Distributed classification: maximum error between the local solution estimates and the optimal solution.}
  \label{fig:logistic_solution}
\end{figure}

The simulation results are reported in Figures~\ref{fig:logistic_cost} and~\ref{fig:logistic_solution}.
It can be seen that, for both the distributed algorithms, the solution and cost error go to zero
(although with different rates).


%

%

\subsection{Distributed Classification via Support Vector Machines}
In this example, we consider again the distributed classification problem and we apply a different strategy to solve it.
\subsubsection{Problem formulation}
If the data are known to be divided into two clusters that can be exactly separated by a hyperplane,
the previous distributed classification problem can be recast as a \emph{hard-margin} SVM problem, i.e.,
\begin{equation}\label{pb:svm}
  \begin{aligned}
    &\m_{w,b} 
    & & \frac{1}{2}\|w\|^2\\
    &\st
    & & \ell_{i,j}(w^\top p_{i,j} + b)\geq 1,&\forall j,~\forall i.
  \end{aligned}
\end{equation}
This is a common-cost problem of the type~\eqref{pb:common_cost}, in which the constraint set of agent $i$ is given by
\begin{equation*}
  X_i = \{(w,b)\mid \ell_{i,j}(w^\top p_{i,j} + b)\geq 1,\, j=1,\dots,m_i\}.
\end{equation*}

\subsubsection{Simulation results}
\begin{figure}[htpb]\centering
  \includegraphics[scale=1]{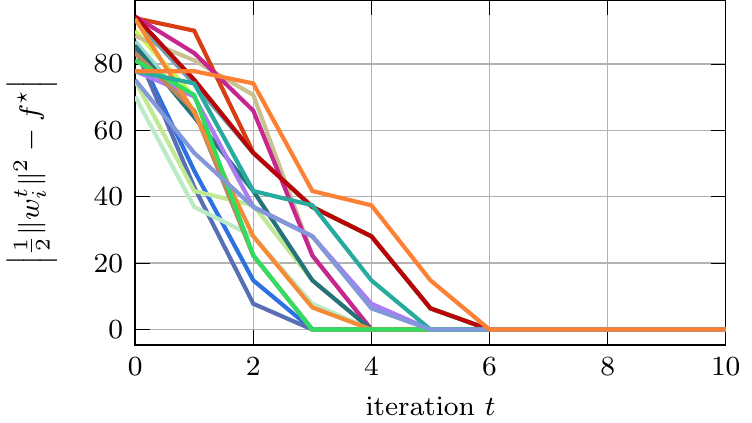}
  \caption{Distributed classification via SVM: cost error between the locally computed solution estimates and the optimal cost.}
  \label{fig:svm_cost}
\end{figure}

\begin{figure}[htpb]\centering
  \includegraphics[scale=1]{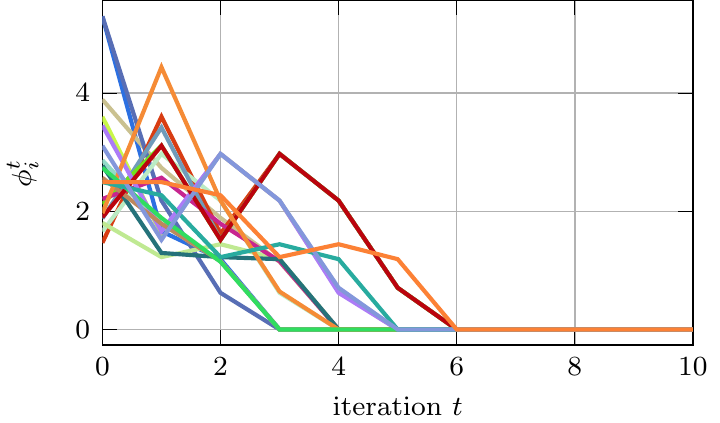}
  \caption{Distributed classification via SVM: maximum constraints violation of the locally computed solution estimates.}
  \label{fig:svm_violation}
\end{figure}

The problem data are generated as described in the previous example
(see Section~\ref{sec:log_loss}).
We apply the Constraints Consensus algorithm, with $N=20$ agents
for $10$ iterations. Figure~\ref{fig:svm_cost} depicts the evolution of the
cost error computed at the local solution estimate of each agent. As
expected from the theoretical analysis (see, e.g., \cite[Chapter 4]{notarstefano2019distributed}),
all the agents converge to a cost-optimal
solution in finite time. Moreover, let us define the maximum
constraint violation of the solution computed by agent $i$ at iteration $t$,
\begin{align*}
  \phi_i^t =1-  \max_{k,j} \ell_{k,j}\big( (w_i^t)^\top p_{k,j} + b_i^t \big).
\end{align*}
In Figure~\ref{fig:svm_violation}, it is shown that $\phi_i^t$ goes to $0$
for all agents, meaning that the solution retrieved by the agents concurrently
satisfies all the constraints of problem~\eqref{pb:svm} in a finite numbe of
iterations.

\subsection{Microgrid Control}

In a microgrid control problem, a network of dynamical systems cooperate in order to compute an optimal
``power profile'' while satisfying both local constraints (e.g., local dynamics or bounds on the state variables)
and global constraints (e.g., compliance with a common resource budget).
We assume that agents are interested in optimizing their profiles over a given discrete-time 
horizon $\{0,1,\ldots, S\}$, for some $S\in\natural$.

\subsubsection{Problem formulation}
Formally, we consider $N$ dynamical units in which each agent $i$ has state $x_i(k) \in \R$ and 
input $u_i(k) \in \R$ for all $k \in \{0,1,\ldots, S\}$. For all $i\in\until{N}$, the states and the
inputs must satisfy a linear dynamics
\begin{align*}
  x_i(k+1) = A_i x_i(k) + B_i u_i(k), \qquad k \in \{0,1,\ldots, S\}
\end{align*}
for given $x_i(0) \in \R$ and given matrices $A_i$ and $B_i$ of suitable dimensions.
The constraint-coupled optimization problem can be cast as follows
\begin{align*}
    \begin{aligned}
        &\m_{x,u}
        & & \sum_{k=1}^S \sum_{i=1}^N \ell_i(x_i(k), u_i(k)) 
        \\
        & \st
        & & x_i(k+1) = A_i x_i(k) + B_i u_i(k), \quad \forall \, k, \: \forall \, i  
        \\
        &
        & & \sum_{i=1}^N \left ( C_i x_i(k) + D_i u_i(k) \right) \le h_k, \quad \forall \, k.
    \end{aligned}
\end{align*}
where $h_k$ are entries of a given vector $h \in \R^S$, $C_i$ and $D_i$ are matrices of suitable dimensions,
and $x$ and $u$ are the collections of all the states and inputs of the agents.
The last line can interpreted as a constraint on the output map of the local systems.

More details can be found in~\cite{notarstefano2019distributed}.

\subsubsection{Simulation results}
We run a comparative study of distributed dual subgradient
and distributed primal decomposition with $N = 20$ agents,
$S = 8$ coupling constraints and $20,000$ iterations.
For both algorithms, we use a diminishing step-size rule $\alpha^t = (1/t)^{0.6}$.

\begin{figure}[htpb]\centering
  \includegraphics[scale=1]{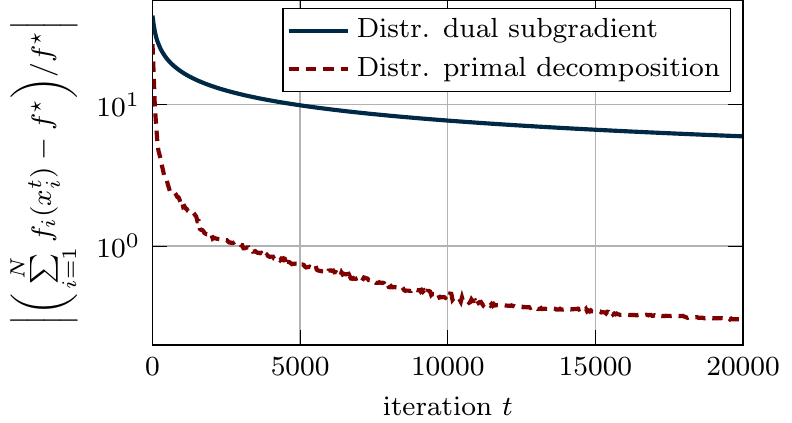}
  \caption{Microgrid control: normalized cost error between the locally computed
    solution estimates and the optimal cost.}
  \label{fig:microgrid_cost}
\end{figure}
\begin{figure}[htpb]\centering
  \includegraphics[scale=1]{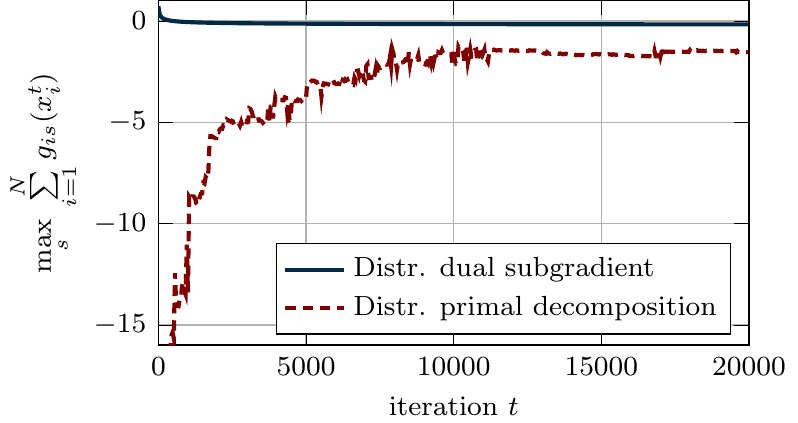}
  \caption{Microgrid control: coupling constraint value for the computed solution estimates.
    Solutions are feasible if the line is below zero.}
  \label{fig:microgrid_coupling}
\end{figure}

The simulation results are reported in Figures~\ref{fig:microgrid_cost} and~\ref{fig:microgrid_coupling}.
For both the distributed algorithms, the cost error goes to zero with sublinear rate.
In this example, for the dual subgradient algorithm, the local solution estimates are
within the coupling constraints after less than $5,000$ iterations, while for the
distributed primal decomposition the solution estimates are always feasible.

\section{Conclusions}
In this paper, we introduced \disropt/, a Python package for distributed optimization
over peer-to-peer networks. The package allows users to define and solve optimization
problems through distributed optimization algorithms.
We presented the software architecture of the package together with simulation results
over example application scenarios.

\bibliographystyle{IEEEtran}
\bibliography{biblio}

\end{document}